\documentstyle{amltd}
\begin{document}
\annalsline{157}{2003}
\received{September 10, 2001}
\startingpage{289}
\def\bye{\end{document}}
 \font\tenrm=cmr10
\input amssym.def
\input amssym.tex
\catcode`\@=11
\font\twelvemsb=msbm10 scaled 1100
\font\tenmsb=msbm10
\font\ninemsb=msbm10 scaled 800
\newfam\msbfam
\textfont\msbfam=\twelvemsb  \scriptfont\msbfam=\ninemsb
  \scriptscriptfont\msbfam=\ninemsb
\def\msb@{\hexnumber@\msbfam}
\def\Bbb{\relax\ifmmode\let\next\Bbb@\else
 \def\next{\errmessage{Use \string\Bbb\space only in math
mode}}\fi\next}
\def\Bbb@#1{{\Bbb@@{#1}}}
\def\Bbb@@#1{\fam\msbfam#1}
\catcode`\@=12

 \catcode`\@=11
\font\twelveeuf=eufm10 scaled 1100
\font\teneuf=eufm10
\font\nineeuf=eufm7 scaled 1100
\newfam\euffam
\textfont\euffam=\twelveeuf  \scriptfont\euffam=\teneuf
  \scriptscriptfont\euffam=\nineeuf
\def\euf@{\hexnumber@\euffam}
\def\frak{\relax\ifmmode\let\next\frak@\else
 \def\next{\errmessage{Use \string\frak\space only in math
mode}}\fi\next}
\def\frak@#1{{\frak@@{#1}}}
\def\frak@@#1{\fam\euffam#1}
\catcode`\@=12
\def\Hol{{\rm Hol}}
\def\ran{{\rm ran\,}}
\def\ep{{}{\hfill $\Box$} \vskip 5pt \par}
\def\={\ = \ }
\font\titr=cmr6 

\title{Norm preserving extensions\\ of holomorphic functions\\
from subvarieties of the bidisk}
 \shorttitle{Holomorphic functions}

 \acknowledgements{The first author was partially supported by the National Science Foundation.  The second author
was partially supported by National Science Foundation
grant DMS-0070639.}
  \twoauthors{Jim Agler}{John E. M\raise.5ex\hbox{\tenrm c}Carthy}
 \shortname{Jim Agler and John E. M\raise2pt\hbox{\titr c}Carthy}
 \institutions{ University of California at San Diego, La Jolla, California\\
{\eightpoint {\it E-mail address\/}: agler@math.ucsd.edu}\\ \vglue6pt 
Washington University, St. Louis, Missouri\\
{\eightpoint {\it E-mail address\/}: mccarthy@math.wustl.edu
}}

\section{Introduction}\label{sec:1}

A basic result in the theory of holomorphic functions of
several complex variables is the following special case of
the work of H. Cartan on the sheaf cohomology on Stein
domains (\cite{car51}, or see \cite{gur} or \cite{hor} for more modern
treatments).

\proclaim{Theorem}\label{thm:1.1}
If $V$ is an analytic variety in a domain of holomorphy
$\Omega$ and if $f$ is a holomorphic function on $V${\rm ,} then
there is a holomorphic function $g$ in $\Omega$ such that
$g=f$ on $V$.
\endproclaim

The subject of this paper concerns an add-on to the structure
considered in Theorem \ref{thm:1.1} which arose in the
authors' recent investigations of Nevanlinna-Pick
interpolation on the bidisk.  The definition for a general
pair $(\Omega, V)$ is as follows.

\numbereddemo{Definition}\label{def:1.2}
Let $V$ be an analytic variety in a domain of holomorphy~$\Omega$.  Say  $V$ has the {\it extension property} if
whenever $f$ is a bounded holomorphic function on $V$, there
is a bounded holomorphic function $g$ on $\Omega$ such
that
\begin{equation}\label{eq:1.3}
g|_V =f \quad \mbox{and}\quad \sup\limits_\Omega |g| =
\sup\limits_V |f|. \speqnu{1.3}
\end{equation}
More generally, if ${\rm Hol}^\infty(V)$ denotes the
bounded holomorphic functions on $V$ and $A\subseteq
{\rm Hol}^\infty(V)$, then we say  $V$ has the
$A$-\/{\it extension property} if there is a bounded holomorphic
function $g$ on $\Omega$ such that (\ref{eq:1.3}) holds
whenever $f\in A$.  \enddemo

Before continuing we remark that in Definition \ref{def:1.2} it
is not essential that $V$ be a variety: interpret $f$ to be
holomorphic on $V$ if $f$ has a holomorphic extension to a
neighborhood of $V$.  Also, in this paper we shall restrict our
attention to the case where $\Omega = {\Bbb D}^2$.  The
authors intend to publish their results on more general cases in a
subsequent paper.  Finally, we point out that the notion in
Definition \ref{def:1.2} is different but closely related to
extension problems studied by the group that worked out the
theory of function algebras in the 60's and early 70's (see e.g.\
\cite{rud69} and \cite{ale69}).  We now describe in some
detail how we were led to formulate the notions in
Definition
\ref{def:1.2}.

The classical Nevanlinna-Pick Theory gives an exhaustive
analysis of the following extremal problem on the disk.  For data
$\lambda_1, \dots, \lambda_n \in {\Bbb D}$ and $z_1, \dots,
z_n\break\in {\Bbb C}$, consider 
\begin{equation}\label{eq:1.4}
\rho \= \inf\  \{ \sup\limits_{\lambda \in {\Bbb D}}\
|\varphi(\lambda)| \ : \ 
\varphi:{\Bbb D}\stackrel{{\rm holo}}{\longrightarrow}
{\Bbb C},\ \varphi(\lambda_i)=z_i \} .\speqnu{1.4}
\end{equation}
Functions $\psi$ for which (\ref{eq:1.4}) is attained are referred
to as {\it extremal} and the most important fact in the whole
theory is that there is only one extremal for given data.  Once
this fact is realized it comes as no surprise that there is a finite
algebraic procedure for creating a formula for the extremal in
terms of the data and the critical value $\rho$
(as an eigenvalue problem) an important result,
not only in function theory \cite{gar81}, but in the model theory
for Hilbert space contractions \cite{foi-fra} and in the
mathematical theory of control \cite{helt87}.

Now, let us consider the associated extremal problem on the
bidisk.  For data $\lambda_i = (\lambda^1_i, \lambda^2_i) \in
{\Bbb D}^2, 1 \le i \le n$, and $z_i \in {\Bbb C}, 1 \le i \le
n$, let
\begin{equation}\label{eq:1.5}
\rho \= \inf \ \{ \sup\limits_{\lambda \in {\Bbb D}^2}\
|\varphi(\lambda)| \ : \
\varphi:{\Bbb D}^2\stackrel{{\rm \ holo}}{\longrightarrow}
{\Bbb C},\ \varphi(\lambda_i)=z_i \}. \speqnu{1.5}
\end{equation}
Unlike   the case of the disk, extremals for (\ref{eq:1.5}) are not
unique.  The authors however have discovered the interesting
fact that there is a polynomial variety in the bidisk on which the
extremals are unique.  Specifically, there exists a polynomial
variety $V_{\lambda, z}\subseteq {\Bbb D}^2$, depending on
the data, and there exists a holomorphic function $f$  defined on
$V_{\lambda, z}$ with the properties that $\lambda_1, \dots,
\lambda_n \in V_{\lambda, z}$:
\begin{equation}\label{eq:1.6}
g|_{V_{\lambda, z}} \= f\quad \mbox{and}\quad
\sup\limits_{{\Bbb D}^2}\ |g| = \sup\limits_{V_{\lambda, z}} \ |f| \speqnu{1.6}
\end{equation}
whenever  $g$  is extremal for (\ref{eq:1.5}).  Furthermore there
is a finite algebraic procedure for calculating $f$ in terms of the
data and the critical value $\rho$ (now, calculating $\rho$ is a
problem in semi-definite programming).
See also the paper \cite{ampa} by Amar and Thomas. 
\advance\theoremcount by 4

Thus, it transpires that there is a unique extremal to
(\ref{eq:1.5}), not defined on all of ${\Bbb D}^2$, but only on
$V_{\lambda, z}$, and that the set of global extremals to
(\ref{eq:1.5}) is obtained by taking the set of norm preserving
extensions of this unique local extremal to the bidisk.  Clearly,
the  nicest possible situation would arise if it were the case that
$V_{\lambda, z}$ had the {\it polynomial extension property} (i.e., Definition~\ref{def:1.2}  holds with $\Omega = {\Bbb
D}^2, V=V_{\lambda, z}$ and $A=$ polynomials), for then the analysis of
(\ref{eq:1.5}) would separate into two independent and
qualitatively different problems: the  analysis of the unique local
extremal, and analysis of the norm preserving extensions from
$V_{\lambda, z}$ to ${\Bbb D}^2$.

Thus, we see that the problem of identifying in the
bidisk, the varieties  that have the $A$-extension property for a given algebra
$A$ arises naturally if one wants to understand the extremal
problem (\ref{eq:1.5}).  It turns out that there is also a purely
operator-theoretic reason to study the $A$-extension property.
This story begins with the famous inequality of von Neumann
\cite{vonN51}.

\proclaim{Theorem}\label{thm:1.7}
If $T$ is a contractive operator on a Hilbert space{\rm ,} then
$$\|p(T)\| \le \sup\limits_{{\Bbb D}}\ |p|$$
whenever $p$ is a polynomial in one variable.
\endproclaim

It was the attempt to explain this theorem that led Sz.-Nagy to
discover his famous dilation theorem \cite{szn53} upon which many of
the pillars of modern operator theory are based (e.g. \cite{szn-foi}).
An extraordinary amount of work has been done
by the operator theory community extending the inequality of
von Neumann, none more elegant than the following result of And\^o
\cite{and63}.

\proclaim{Theorem}\label{thm:1.8}
If $T=(T^1, T^2)$ is a contractive commuting pair of operators on a
Hilbert space{\rm ,} then
\begin{equation}\label{eq:1.9}
\|p(T)\| \le \sup\limits_{{\Bbb D}^2}\ |p| \speqnu{1.9}
\end{equation}
whenever $p$ is a polynomial in two variables.
\endproclaim

We propose in this paper a refinement of Theorem \ref{thm:1.8}
based on replacing (\ref{eq:1.9}) with an estimate
\begin{equation}\label{eq:1.10}
\|f(T)\| \le \sup\limits_{V}\ |f| \speqnu{1.10}
\end{equation}
where $f$ is allowed to be more general than a polynomial and
$V$ is a general subset of the bidisk.  For $V\subseteq
{\Bbb D}^2$, let $\Hol(V)$ denote the functions defined on $V$
that have a holomorphic extension to a neighborhood of $V$ and
let $\Hol^\infty(V)$ denote the set of elements in $\Hol(V)$ that
are bounded on $V$.  Note that if $V$ is a variety and $f$ is
holomorphic on $V$, then a baby theorem in the Cartan theory
would be that there exists a neighborhood $U \supseteq V$ and a
holomorphic function $g$ on $U$  with $f=g|V$, though such a $g$
might well not be  unique.  If we want to form $f(T)$ where $T$ is a
pair of commuting operators, one way would be to define $f(T)$
to be $g(T)$ where $g(T)$  is defined via the Taylor calculus 
\cite{tay70a}.
  Of course, we would need that $\sigma(T) \subseteq V$ (so
that $\sigma(T)\subseteq U$) and, in addition, would want $f(T)$
to depend only on $f$ and not on the particular extension $g$.
This motivates the following definition which makes sense for
{\it arbitrary} sets $V$.
\advance\theoremcount by 2

\numbereddemo{Definition}\label{def:1.11}
If $V\subseteq {\Bbb C}^2$ and $T$ is a commuting  pair of
operators on a Hilbert space, say  $T$ is {\it subordinate to $V$} if
$\sigma(T) \subseteq V$ and $g(T) =0$  whenever $g$ is
holomorphic on a neighborhood of $V$ and $g|V=0$.  If $f\in
\Hol(V)$ and $T$ is subordinate to $V$ {\it define} $f(T)$ by
setting $f(T) =g(T)$ where $g$ is any holomorphic extension of
$f$ to a neighborhood of $V$.  If $T$ is subordinate to $V$ and
(\ref{eq:1.10}) holds for all $f\in \Hol^\infty(V)$, then we say that
  $V$ is a {\it spectral set for $T$}.  More  generally, if $T$ is
subordinate to $V, A \subseteq \Hol^\infty(V)$ and (\ref{eq:1.10})
holds for all  $f\in A$, then we say that   $V$ is an $A$-\/{\it spectral
set for $T$}.
\enddemo

Armed with Definition \ref{def:1.11} it is easy to see that, modulo
some simple approximations, And\^o's theorem is equivalent
to the assertion that ${\Bbb D}^2$ is a spectral set for any pair
of commuting contractions with $\sigma(T)\subseteq
{\Bbb D}^2$.  Thus the following definition seems worthy of
contemplation.

\vglue4pt {\it Definition} 1.12.  
Fix $V\subseteq {\Bbb C}^2$ and let $A\subseteq
\Hol^\infty(V)$.  Say that  $V$ is an $A$-\/{\it von Neumann set} if
$V$ is an $A$-spectral set for $T$ whenever $T$ is a commuting
pair of contractions subordinate to $V$.
\vglue6pt 

We have introduced two properties that a set
$V\subseteq {\Bbb D}^2$ might have relative to a specified
subset $A\subseteq \Hol^\infty(V):V$ might have the
$A$-extension property as in Definition \ref{def:1.2}; or, it might
be an $A$-von Neumann set as in Definition 1.12.
Furthermore, we have indicated the naturalness of these
properties from the appropriate perspectives.  In this paper we
shall show these two notions are actually the same.
Specifically, we have the following result.

\vglue6pt {\elevensc Theorem 1.13}.
{\it Let $V\subseteq {\Bbb D}^2$ and let $A\subseteq
\Hol^\infty(V)$.  Now{\rm ,}  $V$ has the $A$\/{\rm -}\/extension property if and only if
$V$ is an $A$\/{\rm -}\/von Neumann set.}
\vglue6pt
\advance\theoremcount by 2

This theorem will be proved in Section \ref{sec:2} of this paper.

Theorem 1.13 provides a powerful set of tools for
investigating the extension property, namely the techniques of
operator dilation theory.   Specifically, in Section \ref{sec:3} of
this paper we shall show 
that if $V$ is
polynomially convex and $V$ is an $A$-spectral set for a commuting 
pair of $2\times
2$ matrices, then the induced contractive algebra homomorphism of
$\Hol^\infty(V)$ is in fact completely contractive (Proposition 
\ref{prop:3.1}).  It
will then follow via Arveson's dilation theorem \cite{arv69}, operator
model theory, and concrete $H^2$-arguments that $V$ must
satisfy a purely geometric property: $V$ must be balanced.

To define this notion of balanced, we first recall
for the convenience of the
reader some simple notions from the theory of
complex metrics (see  \cite{kra} for an excellent discussion).  
\begin{equation}  \label{eq:1.14}\speqnu{1.14} 
C_U(\lambda_1, \lambda_2) \= \sup \{ d(F(\lambda_1),
F(\lambda_2))\, : \, F:U\to
{\Bbb D},\ F\, {\rm is\,
holomorphic} \} \enspace\
\end{equation}
and 
\begin{equation} \label{eq:1.15}\speqnu{1.15} 
K_U(\lambda_1, \lambda_2) \=
\inf\, \{d(\mu_1, \mu_2) \, : \, {\varphi:{\Bbb D}\to
U},\ {\varphi(\mu_i)=\lambda_i},\
{\varphi\, {\rm is\, holomorphic} } \}  
\end{equation}
where $d(\mu_1, \mu_2) = \left|\frac{\mu_1-\mu_2}{1-\overline
\mu_1 \mu_2}\right|$ is the pseudo-hyperbolic metric on the
disk.  Here, (\ref{eq:1.14}) is referred to as the {\it Carath{\rm \'{\it e}}odory
extremal problem} and  functions for which the supremum in
(\ref{eq:1.14}) is attained are referred to as {\it Carath{\rm \'{\it e}}odory
extremals}.  Furthermore, $C_U(\lambda_1, \lambda_2)$ is
always a metric, the {\it Carath{\rm \'{\it e}}odory metric}.  Likewise,
(\ref{eq:1.15}) is the {\it Kobayashi extremal problem} and
functions for which the infimum is attained are {\it Kobayashi
extremals}.  However, $K(\lambda_1, \lambda_2)$ is in general
not a metric though the beautiful theorem of Lempert \cite{lem81}
asserts that if $U$ is convex, then in fact $K_U$ is a metric,
and
indeed $K_U = C_U$.  In the simple case when $U={\Bbb D}^2$
both (\ref{eq:1.14}) and (\ref{eq:1.15}) are easily solved to yield
the formulas\footnote{We shall use superscripts to denote coordinates in ${\scriptstyle\Bbb D}^2$,
subscripts to distinguish points, or to denote coordinates of a
vector.}
\begin{equation}\label{eq:1.16}
C_{{\Bbb D}^2} (\lambda_1, \lambda_2) =K_{{\Bbb D}^2}
(\lambda_1, \lambda_2) =\max\{d^1, d^2\} \speqnu{1.16}
\end{equation}
where $d^1 =d(\lambda_1^1, \lambda_2^1)$ and
$d^2=d(\lambda_1^2, \lambda_2^2)$.

The formulas (\ref{eq:1.16})  allow one to see that the
description of the extremal functions for (\ref{eq:1.14}) and
({\ref{eq:1.15}) in the case when $U={\Bbb D}^2$ splits naturally
into three cases according as $d^1 > d^2, d_1 = d_2$, or
$d^1 <d^2$.  If  $d^1 >d^2$, then the extremal function for
(\ref{eq:1.14}) is unique: $F(\lambda) = \lambda^1$.  However
when $d^1 >d^2$, there is not a unique extremal for
(\ref{eq:1.15}): any function $\varphi(z) = (z, f(z))$ where
$f:{\Bbb D} \to {\Bbb D}$ satisfies $f(\lambda^1_i) =
\lambda^2_i$ will do.  Likewise when $d^1 <d^2$ the Carath\'eodory
extremal is the unique function $F(\lambda) = \lambda^2$ and any
function $\varphi(z) = (f(z), z)$ where $f:{\Bbb D}\to
{\Bbb D}$ solves $f(\lambda^2_2)= \lambda^1_2$ is a
Kobayashi extremal.  Thus, when $d^1 \not= d^2$, the
Carath\'eodory extremal is unique and the Kobayashi is not.   When
$d^1=d^2$, the reverse is true, the Carath\'eodory extremal is not
unique and the Kobayashi extremal is: either $F(\lambda)
=\lambda^1$ or $F(\lambda) = \lambda^2$ is extremal for
(\ref{eq:1.14}) while $\varphi(z) = (z, f(z))$, where $f$ is the
unique M\" obius mapping  $f:{\Bbb D} \to {\Bbb D}$
satisfying $f(\lambda^1_i) = \lambda^2_i$, is the unique extremal
for (\ref{eq:1.15}).  These considerations prompt the following
definition. \advance\theoremcount by 3

\numbereddemo{Definition}\label{def:1.17}
If $\lambda=(\lambda_1, \lambda_2) = ((\lambda_1^1,
\lambda_1^2),
(\lambda_2^1, \lambda_2^2))$ is a pair of points in
${\Bbb D}^2$, say
$\lambda$ is a {\it balanced pair} if $d(\lambda_1^1,
\lambda_2^1) = d(\lambda_1^2, \lambda_2^2)$.
\enddemo
 
  Thus, the
Kobayashi extremal for a pair of points $\lambda $ is unique if
and only if
$\lambda$ is a balanced pair. Now, if $\lambda$ is pair of points in 
${\Bbb D}^2$,
and
$\varphi$ is extremal for the Kobayashi problem, it is easy to
check that $D=\ran \varphi$ is a totally geodesic one dimensional
complex submanifold of ${\Bbb D}^2$.  Conversely, if $D=\ran
\varphi$ is an analytic disk in ${\Bbb D}^2$ and $D$ is totally
geodesic, then  $\varphi$ is a Kobayashi extremal for any pair of
points in $D$.  Thus, we may assert, based on the observation
following Definition \ref{def:1.17}, that there exists a unique
totally geodesic disk $D_\lambda$ passing through a pair of
points $\lambda$ in ${\Bbb D}^2$ if and only if $\lambda$ is a
balanced \pagebreak pair.  Concretely, it is the set
$$D_\lambda = \{(z, f(z)): z\in {\Bbb D}\}$$
where $f:{\Bbb D}\to {\Bbb D}$ is the unique mapping
satisfying
$f(\lambda^1_i) = \lambda^2_i$.

We now are able to give the promised definition of a balanced
subset of~${\Bbb D}^2$.

\numbereddemo{Definition}\label{def:1.18}
If $V \subseteq {\Bbb D}^2$, say {\it $V$ is balanced} if
$D_\lambda \subseteq V$ whenever $\lambda$ is a balanced pair
of points in $V$.
\enddemo

Note that if $D$ is either an analytic disk or a totally geodesic
disk in ${\Bbb D}^2$, then $D$ is balanced in the sense of
Definition \ref{def:1.18} if and only if $D=D_\lambda$ for some
balanced pair $\lambda$.  For this reason we refer to $D_\lambda$
as  the {\it balanced disk passing through $\lambda_1$ and
$\lambda_2$}.  Note that if $D$ is a balanced disk, then every pair
of  points in $V$ is balanced and also $D=D_\lambda$ for each pair
of points $\lambda \in D \times D$.
The significance of balanced sets in the context of the extension
property on the bidisk will be revealed in Section \ref{sec:3} 
where we shall exploit Theorem 1.13 to give an
operator-theoretic proof of the following result.

\proclaim{Theorem}\label{thm:1.19}
Let $V\subseteq {\Bbb D}^2$ and assume that $V$ is relatively
polynomially convex {\rm (}\/i.e.\ $V^\wedge \cap {\Bbb D}^2 =V$
where $V^\wedge$ denotes the polynomially convex hull of $V${\rm ).}
If $V$ has the polynomial extension property{\rm ,} then $V$ is
balanced.
\endproclaim

It turns out that the property of being balanced is much more
rigid than one might initially suspect.  In Section \ref{sec:4} we
shall investigate this phenomenon by establishing several
geometric properties of balanced sets.  Finally, in Section~\ref{sec:5} of this paper we shall combine this geometric
rigidity of balanced sets with the elementary observation that subsets
of the bidisk with an extension property must be
$H^\infty$-varieties to obtain the following result which gives a
complete classification of the subsets $V$ of the bidisk with the
polynomial extension property (at least in the case when $V$ is
relatively polynomially convex).

\proclaim{Theorem}\label{thm:1.20}
Let $V$ be a nonempty relatively polynomially convex subset of
${\Bbb D}^2$.  $V$ has the polynomial extension property if
and only if $V$ has one of the following forms.
\begin{itemize}
\item[{\rm (i)}]  $V=\{\lambda\}$ for some $\lambda \in {\Bbb D}^2.$

\item[{\rm (ii)}]  $V= {\Bbb D}^2.$

\item[{\rm (iii)}]  $V=\{(z, f(z)|z\in {\Bbb D}\}$ for some holomorphic
$f: {\Bbb D}\to {\Bbb D}.$

\item[{\rm (iv)}]   $V=\{f(z), z)|z\in {\Bbb D}\}$ for some
holomorphic
$f: {\Bbb D}\to {\Bbb D}.$
\end{itemize}

\endproclaim

After this paper was submitted, Pascal Thomas devised an
elegant function theoretic proof of Theorem~\ref{thm:1.19}.
We include his proof in an appendix at the end of the paper.

\section{The equivalence of the von Neumann inequality\\ and
the extension property}\label{sec:2} {\ }

In this section we shall prove Theorem 1.13 from
the introduction.  Accordingly, fix a set
$V\subseteq {\Bbb D}^2$ and a set
$A\subseteq\Hol^\infty(V)$.

One side of Theorem 1.13 is straightforward.  Thus,
assume that $V$ has the $A$-extension property and fix a
commuting pair of contractions $T$ such that $T$ is subordinate to 
$V$.  If $f\in A$
and
$g\in H^\infty({\Bbb D}^2)$ with
$g|V=f$ and $\|g\|_{{\Bbb D}^2} = \|f\|_V$, then
$$\|f(T)\| = \|g(T)\| \le \|g\|_{{\Bbb D}^2} = \|f\|_V.$$
Hence, since $f$ was arbitrarily chosen, $V$ is an $A$-spectral
set for $T$.  Hence since $T$ was arbitrarily chosen, $V$ is an
$A$-von Neumann set. 

The reverse direction of Theorem 1.13 is much more
subtle and will rely on some basic facts about Nevanlinna-Pick
interpolation on the bidisk.  For $n$ distinct points in the bidisk 
$\lambda_1, \dots ,
\lambda_n$,
  let
${\cal  K}_\lambda$ denote the set of $n\times n$ strictly
positive definite matrices, $[k_i (j)]^n_{i, j=1}$, such that $k_i(i)
=1$ for each $i$,
\begin{equation}\label{eq:2.1}
\left[\left(1 - \lambda^{\overline 1}_i \lambda^1_j\right)
k_i(j)\right] \ge 0 \quad \mbox{and}\qquad
\left[\left(1 - \lambda^{\overline 2}_i \lambda^2_j\right)
k_i(j)\right] \ge 0.
\end{equation}
For a proof of the next result see the papers \cite{colwer94}, \cite{baltre98}
or \cite{agmc_bid}, or the book \cite{ampi}.

\advance\theoremcount by 1
\proclaim{Theorem}\label{thm:2.2} \hskip-8pt 
If $\lambda_1, \dots, \lambda_n$ are $n$ distinct points in
${\Bbb D}^2$ and $z_1, \dots, z_n \in {\Bbb C}$ then there
exists $\varphi \in H^\infty({\Bbb D}^2)$ with $\|\varphi
\|_{{\Bbb D}^2} \le 1$ and $\varphi(\lambda_ i) =z_i$ for each
$i$ if and only if
\begin{equation}\label{eq:2.3}
\Big[\left(1 - \overline z_i z_j\right)
k_i(j)\Big] \ge 0, \qquad \mbox{for all}\ \ k\in
{\cal  K}_\lambda. \speqnu{2.3}
\end{equation}
\endproclaim

Theorem \ref{thm:2.2} allows us via a simple argument to
dualize the extremal problem (\ref{eq:1.5}) in the following form.
\advance\theoremcount by 1

\proclaim{Theorem}\label{thm:2.4} \hskip-8pt 
If distinct points $\lambda_1, \dots, \lambda_n\in {\Bbb D}^2
$ and points  $z_1, \dots, z_n \in {\Bbb C}$ are given and
$\rho$ is as in {\rm (\ref{eq:1.5}),}  then
$$\rho = \inf \left\{\sigma\Big|\Big[(\sigma^2 - \overline z_i z_j) k_i
(j)\Big] \ge 0 \quad \mbox{for all}\ k \in {\cal  K}_\lambda\right\}.$$
\endproclaim

Finally, Theorem \ref{thm:2.4} yields the following result which
will provide the key ingredient for completing the proof of
Theorem 1.13.
 
\proclaim{Lemma}\label{lem:2.5}
If $\rho$ is as in {\rm (\ref{eq:1.5})} and $\psi$ is extremal for
{\rm (\ref{eq:1.5}),} then there exist a commuting pair of contractions
$T= (T_1, T_2)$ subordinate to $\{\lambda_1, \dots,
\lambda_n\}$ and such that $\|\psi (T)\| =\rho$. \pagebreak
\endproclaim

\demo{Proof}  We first claim that ${\cal  K}_\lambda$
is compact as a subset of the self-adjoint $n\times n$ complex
matrices equipped with the matrix norm.  To see this we show
that ${\cal  K}_\lambda$ is both bounded and closed.  That
${\cal  K}_\lambda$ is bounded follows  when the
normalization condition $k_i(i) =1$ implies that if $ k\in
{\cal  K}_\lambda$, then
$$\|k\| \le  {\rm tr} \ K = \sum k_i (i) = n.$$

To see that ${\cal  K}_\lambda$ is closed, we argue by
contradiction.  Thus, assume that $\{k^\ell\}$ is a sequence in
${\cal  K}_\lambda, k^\ell \to k$ as $\ell \to \infty$, and
$k\not\in {\cal  K}_\lambda$.  By continuity, $k_i(i) =1$ for
each $i$.  Also by continuity, condition (\ref{eq:2.1}) holds.
Hence since $k\not\in {\cal  K}_\lambda$ it must be the case
that $k$ is not strictly positive definite.  Choose a vector
$v=(v_i)$ with $k v =0$ and $v\not= 0$.  Letting $\Lambda^1$ and
$\Lambda^2$ denote the diagonal matrices whose $(i,i)^{\rm th}$
entries are $\lambda^1_i$ and $\lambda^2_i$ we deduce from
(\ref{eq:2.1}), that if $r=1$ or $r=2$, then,
\begin{eqnarray*}
0&\le& \sum\limits_{i, j} (1-\lambda^{\overline r}_i
\lambda^{ r}_j) k_i (j) v_j \overline v_i\\
&=& \sum\limits_{i, j} k_i(j) v_j \overline v_i - \sum\limits_{i,
j}  k_i(j) (\lambda^r_j v_j)\overline{(\lambda^r_i v_i)}\\
&=& <kv, v> - < k\Lambda^rv, \Lambda^rv>.
\end{eqnarray*}
Now, $kv=0$ and $k$, by continuity, is positive semidefinite.
Hence both $\Lambda^1v$ and $\Lambda^2v$ are in the kernel of
$k$.  Continuing, we deduce by induction that if $m= (m_1, m_2)$
is a multi-index, then $\Lambda^mv = (\Lambda^1)^{m_1}
(\Lambda^2)^{m_2}v$ is in the kernel of $k$.  Finally, if $p$ is
any polynomial in two variables we deduce that
\begin{equation}\label{eq:2.6}
k\ p(\Lambda) v=0. \speqnu{2.6}
\end{equation}
Now $v\not= 0$ so that there exists $i$ such that $v \not= 0$.  On the
other hand $p(\Lambda)$ is the diagonal operator whose
$j-j^{\rm th}$ entry is $p(\lambda_j)$ and $\lambda_1, \dots
\lambda_n$ are assumed distinct so that there is a polynomial
$p$ such that $p(\lambda_i) =1$ and $p(\lambda_j) =0$ for
$j\not= i$.  Hence from (\ref{eq:2.6}) we see that $p(\Lambda) v=
v_i e_i \in \ker k$ which contradicts the fact that $k_i(i) \not=
0$.  This contradiction establishes that ${\cal  K}_\lambda$ is
closed and completes the proof that ${\cal  K}_\lambda$ is
compact.

As an immediate consequence of the compactness of
${\cal  K}_\lambda$ and Theorem \ref{thm:2.4} there exists
$k\in  {\cal  K}_\lambda$ such that
\begin{equation}\label{eq:2.7}
\Big[(\rho^2 - \overline z_i z_j) k_i (j)\Big] \ge 0 \qquad
\mbox{and} \speqnu{2.7}
\end{equation} 
\begin{equation}\label{eq:2.8}
\exists\ w\in {\Bbb C}^n\ \mbox{with}\
w\not= 0\ \mbox{and}\ \sum\limits_{i, j}(\rho^2 -
\overline z_i z_j) k_i (j) \bar w_i\ w_j
= 0. \speqnu{2.8}
\end{equation}

To define $T$ we first define a pair $X=(X^1, X^2)$.  Choose
vectors $k_1, \dots, k_n\break \in {\Bbb C}^n$ such that $<k_i, k_j>\
= k_i(j)$ for all $i$ and $j$.  Since $k$ is strictly positive
definite, the formulas
$$X^r k_i = \overline{\lambda^r_i} k_i\qquad 1\le i \le n, \ r=1,
2$$
uniquely define a commuting pair of $n\times n$ matrices
$X=(X^1, X^2)$.  Set $T=(X^{1*}, X^{2*})$.  Since $X^1$ and $X^2$
share a set of $n$ eigenvectors with corresponding eigenvalues
$\overline \lambda_1, \dots \overline \lambda_n$ it is clear
that $X$ is subordinate to $\{\overline \lambda_1, \dots,
\overline \lambda_n\}$.  Hence $T$ is subordinate to $\{
\lambda_1, \dots,  \lambda_n\}$.  Noting that
(\ref{eq:2.1}) implies that both $X^1$ and $X^2$ are contracting
we see that $T$ is a contractive pair.  Finally, note that if $\psi$
is extremal for (\ref{eq:1.5}) and $\psi^{\breve{}}$ is defined by
$\psi^{\breve{}}(\lambda) = \overline{ \psi(\overline \lambda)}$, then
$\phantom{\sum^{\raise2pt\hbox{$1$}}}$
$$\psi^{\breve{}}(X) k_i = \psi^{\breve{}}(\overline \lambda_i) k_i =
\overline{\psi(\lambda_i)} k_i = \overline z_i k_i.$$
Hence (\ref{eq:2.7}) implies that $\|\psi^{\breve{}}(X) \| \le \rho$ and
(\ref{eq:2.8}) implies that $\|\psi^{\breve{}}(X)\| \ge \rho$.  Hence
$\|\psi^{\breve{}}(X)\|=\rho$.  But $\psi^{\breve{}}(X) 
= \psi(T)^*$ so that 
$\|\psi(T)\|=\rho$.  This establishes Lemma \ref{lem:2.5}.
\enddemo

We now are ready to complete the proof of Theorem 1.13.  Thus,
assume that $V$ is an $A$-von Neumann set and fix $f\in
A$.  We need to show that there exists $g\in H^\infty({\Bbb D}^2)$
with $g|_V=f|_V$ and $ \|g\|_{{\Bbb D}^2} = \|f\|_V$.

Choose a dense sequence $\{\lambda_i\}^\infty_{i=1}$ in $V$.
For each $n\ge 1$ consider the extremal problem
\begin{equation}\label{eq:2.9}
\rho_n \= \inf\limits_{
{\varphi:{\Bbb D}^2\to
{\Bbb D}} }
\{ \|\varphi\|_{{\Bbb D}^2} \, : \, 
{\varphi(\lambda_i) =
f(\lambda_i)\, {{\rm for}}\, i\le n},\,
{\varphi\, {\rm is\, holomorphic} } \}. \hskip.25in \speqnu{2.9}
\end{equation}
 
If $\psi_n$ is chosen extremal for (\ref{eq:2.9}), then
\begin{eqnarray*}
\|\psi_n\|_{{\Bbb D}^2} &\stackrel{{\rm (i)}}{=}& \rho_n\\
&\stackrel{{\rm (ii)}}{=}& \|\psi_n(T_n)\|\\
&\stackrel{{\rm (iii)}}{=}& \|f(T_n)\|\\
&\stackrel{{\rm (iv)}}{\leq}& \|f\|_V.
\end{eqnarray*}
Here, (i) holds since $\psi_n$ is extremal for (\ref{eq:2.9}), $T_n$
is the commuting pair of contractions subordinate to
$\{\lambda_1, \dots, \lambda_n\}$ whose existence is
guaranteed by Lemma~\ref{lem:2.5}, (ii) holds by Lemma~\ref{lem:2.5}, (iii) holds since $T_n$ is subordinate to
$\{\lambda_1, \dots, \lambda_n\}$ and $\psi_n(\lambda_i)
=f(\lambda_i)$ for $i\le n$, and (iv) holds from the assumption
that $V$ is an $A$-von Neumann set and the fact that $T_n$ is a
contractive pair\break ($T_n$ is subordinate to $V$ since $T_n$ is
subordinate to $\{\lambda_1, \dots, \lambda_n\} \subseteq V$.)

Summarizing, in the previous paragraph we have shown that for
each $n\ge 1$, there exists $\psi_n\in H^\infty({\Bbb D}^2)$
with
\begin{equation}\label{eq:2.10}
\|\psi_n\|_{{\Bbb D}^2} \le \|f\|_V\qquad \mbox{and} \speqnu{2.10}
\end{equation}
\begin{equation}\label{eq:2.11}
\psi_n(\lambda_i) = f(\lambda_i) \quad i \le n. \speqnu{2.11}
\end{equation}
Evidently, either by a uniform family argument or a weak-*
compactness in $H^\infty$ argument, (\ref{eq:2.10}) implies that
there exists $g\in H^\infty({\Bbb D}^2)$ with $\psi_n \to g$
pointwise on ${\Bbb D}^2$ and $\|g\|_{{\Bbb D}^2} \le
\|f\|_V$.  By (\ref{eq:2.11}) we also have that
$g(\lambda_i) = f(\lambda_i)$ for all $i$.  Since
$\{\lambda_i\}^\infty_{i=1}$ was chosen dense in $V$ it follows
that $g|V=f$.  We have shown that $g$ exists with the desired
properties and the proof of Theorem 1.13 is complete.
\hfill\qed

\vglue-6pt
\section{Sets with the polynomial extension property are balanced}\label{sec:3}
 \vglue-6pt

In this section we shall prove Theorem \ref{thm:1.19}  from the
introduction.  
In the statement of our first result, 
note that $\lambda^1$ and $\lambda^2$ are the coordinate
functions, and $\lambda_1$ and $\lambda_2$ are points in ${\Bbb D}^2$.
We use $V^-$ to denote the closure of $V$.

\proclaim{Proposition}\label{prop:3.1}
Let $V \subseteq {\Bbb D}^2${\rm ,} assume that $V$ has the
polynomial extension property{\rm ,} and let $\lambda_1, \lambda_2
\in V$ with $\lambda_1 \not= \lambda_2$.  If $\tau$ is a
contractive\break $2$\/{\rm -}\/dimensional representation of $P(V^-)$ with
$(\tau(\lambda^1), \tau(\lambda^2)) = (\lambda_1,
\lambda_2)${\rm ,} then $\tau$ is completely contractive.
\endproclaim

\demo{Proof}  Let $\Omega$ denote the unit ball of the
complex $n\times n$ matrices and fix an $n\times n$ matrix $p$
of polynomials in two variables with
\begin{equation}\label{eq:3.2}
\|p(\lambda)\| <1\quad \mbox{ for all}\ \ \lambda \in V. \speqnu{3.2}
\end{equation}
Let $T=(\tau(\lambda^1), \tau(\lambda^2))$ be a pair of commuting
contractions on ${\Bbb C}^2$.  Proposition
\ref{prop:3.1} will be established if we can show that
\begin{equation}\label{eq:3.3}
\|p(T)\| \le 1. \speqnu{3.3}
\end{equation}

We claim that

\begin{equation}\label{eq:3.4}
d_\Omega(p(\lambda_1), p(\lambda_2)) \le
d_{{\Bbb D}^2}(\lambda_1, \lambda_2), \speqnu{3.4}
\end{equation}
where $d_\Omega$ (respectively, $d_{{\Bbb D}^2}$)
is the Carath\'eodory metric in $\Omega$ (resp.,
${{\Bbb D}^2}$).
To see (\ref{eq:3.4}), fix $\epsilon >0$.  Choose a polynomial $q$
such that $q:\Omega \to {\Bbb D}$ and
$d_{{\Bbb D}}(q(p(\lambda_1)), q(p(\lambda_2))) >
d_\Omega(p(\lambda_1), p(\lambda_2)) - \epsilon$.  Thus
$q\circ p$ is a  polynomial and (\ref{eq:3.2}) implies that
$\sup\limits_V |q\circ p| \le 1$.  Since $V$ has the polynomial
extension property, there exists $\varphi \in
H^\infty({\Bbb D}^2)$ with $\sup\limits_{{\Bbb D}^2}
|\varphi| \le 1$ and $\varphi(\lambda_i) = q(p(\lambda_i))$.
Hence, for each $\varepsilon >0$,
\begin{eqnarray*}
d_\Omega(p(\lambda_1), p(\lambda_2)) - \epsilon
&\ <  \ & d_{{\Bbb D}}(q(p(\lambda_1)), q(p(\lambda_2))) \\
&\le&
d_{{\Bbb D}^2}(\lambda_1, \lambda_2),
\end{eqnarray*}
which establishes (\ref{eq:3.4}).

To see (\ref{eq:3.3}), choose points $z_1, z_2\in {\Bbb D}$
such that $d(z_1, z_2) = d_\Omega(p(\lambda_1),
p(\lambda_2))$ and let $f:{\Bbb D}\to \Omega$ with $f(z_1) =
p(\lambda_1)$ and $f(z_2) = p(\lambda_2)$.  Evidently,
(\ref{eq:3.4}) implies that there exists $\psi:{\Bbb D}^2\to
{\Bbb D}$ such that $\psi(\lambda_1)=z_1$ and
$\psi(\lambda_2) =z_2$.  Hence, since $T$ is a pair of commuting
contractions $(\tau$ is assumed contractive), it follows
from And\^o's theorem (Theorem \ref{thm:1.8}) that $\psi(T)$ is a
contraction.  Since $\psi(T)$ is a contraction and
$f:{\Bbb D}\to \Omega$, it follows that $f(\psi(T))$ is a
contraction. (This follows from the Sz.-Nagy dilation theorem \cite{szn53},
 as observed by\break Arveson \cite{arv69}: in modern language, the disk
is a complete spectral set for any contraction.)  
But by construction, $f(\psi(T))
= p(T)$.  Thus, (\ref{eq:3.3}) holds and the proof of Proposition
\ref{prop:3.1} is complete.
\enddemo

Our next result exploits the Arveson extension theorem
\cite{arv69},  the
Stinespring representation theorem \cite{sti55}, and the Sarason
interpretation of semi-invariant subspaces \cite{sar65b} to interpret
the completely contractive representation of Proposition
\ref{prop:3.1}.  
We shall say that a subnormal pair $S$ has extension spectrum in
$\Gamma$ if there is a commuting normal pair $N$ whose spectral
measure is supported by $\Gamma$ and such that the restriction of
$N$ to an invariant subspace is unitarily equivalent to $S$.
\advance\theoremcount by 3

\proclaim{Proposition}\label{prop:3.5}
Let $V\subseteq {\Bbb D}^2$ and let $\Gamma$ denote the
Shilov boundary of $P(V^-)$.  Let ${\cal  H}$ be a Hilbert space
and let $\tau:P(V^-) \to {\cal  L}({\cal  H})$ be a
completely contractive representation.  There exists a Hilbert
space ${\cal K}$ and a subnormal pair $S$ with
extension spectrum in $\Gamma$ such that for all $\varphi
\in P(V^-), {\cal  K} \ominus {\cal  H}$ is invariant for
$\varphi(S)$ and $\tau(\varphi)= P_{{\cal  H}}
\varphi(S)|{\cal  H}$.
\endproclaim

{\it Proof}.
By the theorems of Arveson and Stinespring, there is a Hilbert
space
${\cal G}$
containing ${\cal H}$ and a representation
$\pi: C(\Gamma) \to {\cal L}({\cal G})$
such that
$$
\tau(\phi) \= P_{\cal H} \, \pi(\phi) \, |_{\cal H}
\qquad \hbox{for all }\ \phi \, \in \, P(V^-) .
$$
Let $N = (\pi(\lambda^1), \pi(\lambda^2))$. Then the spectrum of
$N$ is contained in $\Gamma$. As $\tau$ is a representation, it
follows from Sarason's lemma that ${\cal H}$ is semi-invariant
for $\pi(P(V^-))$. This means that there exists a superspace
${\cal K} \supseteq {\cal H}$
such that ${\cal K}$ and ${\cal K} \ominus {\cal H}$
are both invariant for $\pi(P(V^-))$. Let $S$ be $N\, |_{\cal K}$.
\hfill\qed\vglue6pt

Armed with propositions we are now ready to commence the
proof of Theorem \ref{thm:1.19}.  Accordingly, fix $V\subseteq
{\Bbb D}^2$ and assume that
\begin{equation}\label{eq:3.6}
V^\wedge \cap {\Bbb D}^2 =V \speqnu{3.6}
\end{equation}
and
\begin{equation}\label{eq:3.7}
V \quad \mbox{has the polynomial extension property.} \speqnu{3.7}
\end{equation}
Fix $d>0$ and a pair of points $\lambda = (\lambda_1,
\lambda_2)$ with $\lambda_i\in V$ and with the property that
$d(\lambda_1^1, \lambda_2^1) =d =d(\lambda_1^2,
\lambda_2^2)$.  Define a  pair of operators $T=(T_1, T_2)$ by
the formulas, $T_rk_i = \overline \lambda^r_i k_i, i=1,2$, and
$r=1, 2$, where $k_i, i=1, 2$, is a pair of unit vectors in
${\Bbb C}^2$ with the property that $|<k_1, k_2>|^2
=1-d^2$.  Noting  that $d=d(\lambda_1, \lambda_2)$, we see
from \cite{ag91} that 
\begin{equation}\label{eq:3.8}
\|\varphi^{\breve{}}(T)\|=1, \speqnu{3.8}
\end{equation}
whenever $\varphi$ is an extremal for the Carath\'eodoty problem for the pair $\lambda$.  (Recall that
$\varphi^{\breve{}}$ is defined by
$\varphi^{\breve{}}(\lambda)=\overline{\varphi(\overline
\lambda)}$.)

Now observe that $\|T_1\| =1$ and $\|T_2\|=1$.  Hence, since 
(\ref{eq:3.7}) holds we
see from Theorem 1.13  that $V$ is a
polynomial von Neumann set. Consequently, if we define
$\tau$ by the formula $\tau(\psi)^* = \psi^{\breve{}}(T)$,
then $\tau$ is a contractive representation of $P(V^-)$.  It
follows from Proposition \ref{prop:3.1} that $\tau$ is
completely contractive.  Consequently, if we let $\Gamma$
denote the Shilov boundary of $P(V^-)$ we obtain from
Proposition \ref{prop:3.5} a commuting subnormal pair of
operators $S$ with extension spectrum in $\Gamma$ and a
$2$-dimensional invariant subspace ${\cal  M}$ for $S$
such that $T_i \cong S^*_i |{\cal  M}$.

Now observe that since $\|T_1\| =1$, there exists a vector
$\gamma \in {\cal  M}$ with $\|\gamma\|=1$ such that
$$P_{{\cal  M}} S_1(S^*_1 \gamma) = T^*_1 T_1 \gamma
=\gamma.$$
Since $\|S_1\| \le 1$, it follows that $S^*_1 \gamma \in
{\cal  M}$.  Thus both $\gamma$ and $S^*_1 \gamma$ are in
${\cal  M}$ and it certainly cannot be the case that they
are linearly dependent, for otherwise $T_1$ would have a
unimodular eigenvalue, contradicting the
fact that
$\sigma(T_1) = \{\lambda^1_1, \lambda^1_2\} \subseteq
{\Bbb D}^2$.  It follows that if we let ${\cal  N}$ denote
the closed linear span of $\{S^{n_1}_1\,
S^{n_2}_2(S^*_1\, \gamma):n_1, n_2 \ge 0\}$ and set $R_i
= S_i|_{\cal  N}$, then $R = (R_1, R_2)$ is a commuting
subnormal pair with extension spectrum in
$\Gamma$.  Furthermore, $R$ is cyclic with cyclic vector
$S^*_1 \gamma, R_1(S^*_1\gamma) =\gamma, {\cal  M}
\subseteq {\cal  N}, {\cal  M}$ is invariant for $R_1^*$
and $R_2^*$, and $T_i= R^*_i|_{\cal  M}$.

The facts in the preceding paragraph imply that there is a
model for $T$ and $R$ of the following type: for some
probability measure $\mu$ with support on $\Gamma$,
${\cal  N} = H^2(\mu)$, the closure of the polynomials in
$L^2(\mu)$; $R_1 = M_{\eta^1}$ and
$R_2=M_{\eta^2}$, multiplication by the coordinate
functions; evaluation at $\lambda_1$ and evaluation at
$\lambda_2$ on polynomials extend by continuity to
continuous linear functionals on $H^2(\mu)$ represented, say,
by the vectors $k_{\lambda_1}$ and $k_{\lambda_2};
{\cal  M} = {\rm span}\{k_{\lambda_1}, k_{\lambda_2}\}$;
$S^*_1\gamma =1; T_r k_{\lambda_1}=\overline\lambda^r_i
k_{\lambda_i}$.  Also, note that (\ref{eq:3.8}) in this model
asserts that
\begin{equation}\label{eq:3.9}
\|M^*_\varphi|{\cal  M}\| =1 \speqnu{3.9}
\end{equation}
whenever $\varphi$ is an extremal for the Carath\'eodory problem for the pair $\lambda$.

Now let $D_\lambda = \{(z, m(z)): z\in {\Bbb D}\}$ be a
parametrization of the balanced disk passing through the
points $\lambda_1$ and $\lambda_2$.  If functions $\varphi_1$
and $\varphi_2$ are defined by the formulas $\varphi_1(\eta)
= \eta^1$ and $\varphi_2(\eta)\, = \, m^{-1}(\eta^2)$, then
$\varphi_1$ and $\varphi_2$ are extremal for the
Carath\'eodory problem for the pair $\lambda$ and have the
additional property that
\begin{equation}\label{eq:3.10}
\varphi_1(\lambda_i) = \lambda^1_i = \varphi_2(\lambda_i)
\ \ \mbox{for}\ \ i=1, 2. \speqnu{3.10}
\end{equation}
Now, since $\varphi_1$ is extremal, (\ref{eq:3.9}) and the fact
that $S^*_1\gamma=1$ imply that
$$
\int|\varphi_1|^2 d \mu = \|M_{\varphi_1}1\|^2 =1.
$$
Hence, \pagebreak $|\varphi_1(\eta) | =1\  \mu$ a.e. But (\ref{eq:3.10})
implies that $M^*_{\varphi_1}|_{\cal  M} =
M^*_{\varphi_2}|_{\cal  M}$ so that $M_{\varphi_2}$ also
attains its norm on 1.  Hence we can also deduce that
$|\varphi_2(\eta)| =1\ \mu$ a.e. Since 
$|\varphi_1|,\ |\varphi_2|,$ and $\left| \frac{\phi_1 +
\phi_2}{2} \right|$ are all $ 1 \ \mu$ a.e., we
have that $\mu $ is actually supported on $\partial D$.

Summarizing, we have shown that if $\lambda =(\lambda_1,
\lambda_2)$ is a balanced pair of points in $V$ and
$D_\lambda$ is the balanced disk passing through the pair,
then there is a probability measure $\mu$ supported on
$\Gamma \cap \partial D_\lambda$ such that $\lambda_1$
and $\lambda_2$ are bounded point evaluations for
$H^2(\mu)$.  Changing variables to ${\Bbb D}$ and invoking
the classical result of Kolmogorov, Szeg\H{o}, and Krein
yields that in fact every point in $D_\lambda$ is a bounded point 
evaluation of $H^2(\mu)$. As bounded point
evaluations are
necessarily in the polynomial convex hull of the support of
$\mu$, we conclude via (\ref{eq:3.6}) that $D_\lambda
\subseteq V$ and the proof of Theorem \ref{thm:1.19} is
complete.
\hfill\qed

\vglue-6pt

\section{Some geometric properties of balanced
sets}\label{sec:4}
\vglue-6pt

In this section we shall derive several simple facts about
balanced subsets of the bidisk.  The main point is that
if there is more than one balanced disk in a balanced set, then
the set must fill in to the entire bidisk.  This will provide us
with a lot of mileage in Section \ref{sec:5} when we prove the
classification theorem (Theorem \ref{thm:1.20}) from Section
\ref{sec:1}.

\proclaim{Lemma}\label{lem:4.1}
Let $B$ be a  balanced subset of ${\Bbb D}^2$.  If $D_1$ is a 
balanced disk in $B$
and $\lambda \in B\setminus D_1${\rm ,} then there exists a balanced disk $D_2$ in
$B$ with $\lambda \in D_2$ and $D_1\cap D_2 \not=
\emptyset$.
\endproclaim

\demo{Proof}  Let $D_1=\{(z, f(z)): z\in {\Bbb D}\}$
be a parametrization of $D_1$ with $f$ an automorphism of
the disk.  If we define a function $\rho$ on ${\Bbb D}$ by the
formula
$$\rho(z) =d(\lambda^1, z)  -  d(\lambda^2, f(z)),$$
then the facts that $\lambda \not\in D_1$ and $f$ is an
isometry guarantee that $\rho(\lambda^1)<~0$ and
$\rho(f^{-1}(\lambda^2)) >0$.  Consequently, the
connectedness of ${\Bbb D}$ implies that there exists a
point $z_0\in {\Bbb D}$ with the property that $\rho(z_0)
=0$.  But since $\rho(z_0) =0$, $((\lambda^1, \lambda^2), (z_0,
f(z_0)))$ is a balanced pair of points in $B$.  Since $B$  is
assumed balanced, it follows that if we let $D_2$ denote the
balanced disk that passes through these points, then $D_2$
lies in $B$.  Finally, the fact that $(z_0, f(z_0)) \in D_1 \cap
D_2$ completes   the proof of Lemma \ref{lem:4.1}.
\enddemo

\proclaim{Lemma}\label{lem:4.2}
For $\tau \in \partial {\Bbb D}$ and $\alpha \in {\Bbb D}${\rm ,}
let $m_{\tau, \alpha}$ denote the automorphism of
${\Bbb D}$ defined by $m_{\tau, \alpha}(z) = \tau\
\frac{z-\alpha}{1-\overline \alpha z}$.  If we view $m_{\tau,
\alpha}$  as a mapping from ${\Bbb D}^-$ into
${\Bbb D}^-${\rm ,} then exactly one of the following four cases
occurs.
$$
\begin{array}{rl}
 {\rm (i)}&\quad  \hbox{$\tau=1$ and $\alpha=0$.}
\\
{\rm (ii)}&\quad \hbox{$\tau\not= 1$ and $|1-\tau| =2|\alpha|$.}
\\
{\rm (iii)}&\quad \hbox{$|1-\tau| <2|\alpha|$.}
\\
{\rm (iv)}&\quad \hbox{$|1-\tau| >2|\alpha|$.}
\end{array}
$$
 Furthermore{\rm ,}
\begin{itemize}
\item[{\rm (i)}] holds if and only if every point in ${\Bbb D}^-$ is a fixed
point for $m_{\tau, \alpha}$.
\item[{\rm (ii)}] holds if and only if  $m_{\tau, \alpha}$ has exactly one
fixed point on $\partial {\Bbb D}^-$.
\item[{\rm (iii)}] holds if and only if  $m_{\tau, \alpha}$ has exactly two
fixed points on $\partial {\Bbb D}^-$.
\item[{\rm (iv)}] holds if and only if  $m_{\tau, \alpha}$ has exactly one
fixed point in ${\Bbb D}$.
\end{itemize}
\endproclaim 

We eschew the proof of Lemma \ref{lem:4.2}.

\proclaim{Lemma}\label{lem:4.3}
Let $B$ be a subset of ${\Bbb D}^2$ with the balanced point
property.  If $D_1$ and $D_2$ are two distinct balanced disks
in $B$ and $\lambda \in D_1 \cap D_2${\rm ,} then there exists a set
$U$ in ${\Bbb C}^2$ such that $U$ is open in ${\Bbb C}^2$
and $\lambda \in U \subseteq B$.
\endproclaim

\demo{Proof}  Let $B$ be a balanced subset of
${\Bbb D}^2$ and let $D_1$ and $D_2$ be as in the lemma.
By composing with an appropriate automorphism of
${\Bbb D}^2$ we can reduce the lemma to the special case
when $(0, 0) \in D_1 \cap D_2, D_1 = \{(z, z):z\in {\Bbb D}\}$,
and for some $\omega \in \partial {\Bbb D}, D_2 = \{(z,
\omega z): z \in {\Bbb D}\}$.  For $\tau \in \partial
{\Bbb D}$ and $\alpha \in {\Bbb D}$, let $m_{\tau,
\alpha}$ denote the automorphism of ${\Bbb D}$ defined by
$m_{\tau, \alpha}(z) = \tau \frac{z-\alpha}{1-\overline \alpha
z}$ and let $D_{\tau, \alpha}$ denote the balanced disk defined
by $D_{\tau, \alpha} =\{(z, m_{\tau, \alpha}(z)): z\in
{\Bbb D}\}$.  Evidently, what we would like to show is that
for each sufficiently small $\lambda \in {\Bbb D}^2$, there
exist $\tau \in \partial {\Bbb D}$ and $\alpha \in
{\Bbb D}\setminus\{0\}$ such that $\lambda \in D_{\tau, \alpha}, D_{\tau,
\alpha} \cap D_1\not= \emptyset$, and $D_{\tau, \alpha} \cap
D_2 \not= \emptyset$.  Equivalently, we want to
show that for each sufficiently small $\lambda \in
{\Bbb D}^2$, there exist $\tau \in \partial {\Bbb D}$ and
$\alpha \in {\Bbb D}\setminus\{0\}$ such that
\begin{equation}\label{eq:4.4}
m_{\tau, \alpha}(\lambda^1) = \lambda^2,  \speqnu{4.4}
\end{equation}
\begin{equation}\label{eq:4.5}
m_{\tau, \alpha}(z) = z \ \ \mbox{for some}\ \ \ z \in
{\Bbb D}, \speqnu{4.5}
\end{equation}
and
\begin{equation}\label{eq:4.6}
m_{\tau, \alpha}(z) = \omega z \ \ \mbox{for some}\ \ \ z \in  
{\Bbb D}. \speqnu{4.6}
\end{equation}

Now observe that (\ref{eq:4.5}) asserts that $m_{\tau, \alpha}$
has a fixed point in ${\Bbb D}$.  Similarly, (\ref{eq:4.6})
asserts that $\overline \omega m_{\tau, \alpha} =
m_{\overline \omega \tau, \alpha}$ has a fixed point in
${\Bbb D}$.  Furthermore, the general  $m_{\tau, \alpha}$
satisfying (\ref{eq:4.4}) is given in terms of a free unimodular
constant, $e^{it}$, by
\begin{equation}\label{eq:4.7}
\tau=\frac{e^{it} - \overline\lambda^1 \lambda^2}{1-e^{it}
\lambda^1 \overline\lambda^2}\quad \mbox{and}\quad
\alpha=\frac{e^{it} \lambda^1- \lambda^2}{e^{it}-
\overline\lambda^1 \lambda^2}. \speqnu{4.7}
\end{equation}
An examination of Case (iv) in Lemma \ref{lem:4.2} thus
reveals that Lemma \ref{lem:4.3} will follow if we can show
that for all sufficiently small $\lambda \in {\Bbb D}^2$,
there  exists $e^{it}$ such that if $\tau$ and $\alpha$ are
defined as in (\ref{eq:4.7}), then $|1-\tau| > 2|\alpha| > 0$ and
$|1-\overline \omega \tau| >2|\alpha| > 0$.  Noting that if
$\lambda$ is small, so also is $\alpha$ and that if $\lambda$ is
small, $\tau$ and $\overline \omega \tau$ are close to $e^{it}$
and $\overline \omega e^{it}$, we see that the lemma follows
by choice of $e^{it}$ so that neither of these latter quantities
is close to $1$, and so that $e^{it} \lambda^1 \neq \lambda^2$.
\enddemo
\advance\theoremcount by 4

\proclaim{Lemma}\label{prop:4.8}
If $B$ is a balanced subset of ${\Bbb D}^2$ and $B$ has
nonempty interior{\rm ,} then $B={\Bbb D}^2$.
\endproclaim
 
{\it Proof}.  Since automorphisms of ${\Bbb D}^2$
map balanced sets to balanced sets, we can assume that $(0,
0)$ is in the interior of $ B$.  
Let $E = \{\lambda \in {\Bbb D}^2\ :\ \
|\lambda^1| = |\lambda^2|\}$.  We claim that
\begin{equation}\label{eq:4.9}
E\subseteq B. \speqnu{4.9}
\end{equation}
To see (\ref{eq:4.9}) fix $\lambda \in E$.  Since $(0, 0)$ is in
the interior of $B$,
there exists $\varepsilon \ >0$ such that $\varepsilon  \lambda \in B$.  Now, $B$
is balanced and $((0, 0), \varepsilon  \lambda)$ is a balanced pair.  Hence
$D=\left\{\left( z, \frac{\lambda^2}{\lambda^1}\ z\right) \ |\ \
z\in {\Bbb D}\right\} \subseteq B$.  In particular, $\lambda
\in B$ and we have shown that (\ref{eq:4.9}) holds.

To see that $B={\Bbb D}^2$, fix $\lambda \in
{\Bbb D}^2\setminus E$, choose an automorphism
$m:{\Bbb D}\to {\Bbb D}$ such that $m(\lambda^1) =
\lambda^2$ and let $D=\{(z, m(z))\ |\ \ z\in {\Bbb D}\}$.
Evidently, since $\lambda \in D$, Lemma \ref{prop:4.8}
will follow if we can show that $D\subseteq B$.  In turn,
$D\subseteq B$ will follow from (\ref{eq:4.9}) and the fact
that $B$ is balanced if we can show that there exist two
distinct points in $D\cap E$.  Now, since $\lambda \not\in E,
m(0) \not=0$.  Hence sometimes $|z| < |m(z)|$ and sometimes
$|z| > |m(z)|$.  It follows from the intermediate value theorem
that there is in fact a continuum of points $z$ such that $|z| =
|m(z)| $, i.e., $(z, m(z)) \in D\cap E$.

This completes the proof of Lemma \ref{prop:4.8}.
\hfill\qed
\advance\theoremcount by 1

\proclaim{Proposition}\label{prop:4.10}
If $B$ is a balanced subset of ${\Bbb D}^2$ and $B$ contains
a balanced pair of points{\rm ,} then either $B$ is a balanced disk or
$B={\Bbb D}^2$.
\endproclaim

{\it Proof}  $B$ contains a balanced pair, say,
$\lambda$.  Since $B$ is balanced, $D_\lambda \subseteq B$.  If
$B=D_\lambda$, then $B$ is a balanced disk.  Thus, the
proposition will follow if it is the case that the
assumption $B\not = D_\lambda$ implies $B={\Bbb D}^2$.

Accordingly, assume that $B\setminus D_\lambda \not=
\emptyset$.  By Lemma \ref{lem:4.1}, there exists a balanced disk
$D$ with $D\not= D_\lambda$ and $D\cap D_\lambda \not=
\emptyset$.  By Lemma \ref{lem:4.3}, $B$ has nonempty interior.
Hence Lemma \ref{prop:4.8} implies that
$B={\Bbb D}^2$ and the proof of Proposition
\ref{prop:4.10} is complete.
\hfill\qed

\section{Which sets have the extension
property?}\label{sec:5}

In this section we shall prove 
Theorem \ref{thm:1.20} from Section
\ref{sec:1} of this paper.

\proclaim{Lemma}\label{lem:5.1}
Let $V\subseteq {\Bbb D}^2$ and assume that $V^\wedge
\cap {\Bbb D}^2 =V$.  If $V$ has the polynomial extension
property{\rm ,} then either $V$ consists of a single point or $V$ has
no isolated points.
\endproclaim

\demo{Proof}  The conclusion of the lemma is logically
equivalent to the assertion that if $V$ has an isolated point
then $V$ is a point.  Accordingly, assume $\lambda_0\in V$ is
an isolated point of $V$.  Now, the maximal ideal space of
$P(V^-)$ is $V^\wedge$ and by hypothesis $V^\wedge \cap
{\Bbb D}^2=V$.  Hence since $\lambda_0$ is an isolated
point in $V, \lambda_0$ is an isolated point in the maximal
ideal space of $P(V^-)$.  It follows from the Shilov idempotent
theorem (see  {e.g.}\ \cite{alewer})
that there exists $e\in P(V^-)$ with
$e(\lambda_0)=1$ and $e(\lambda)=0$ for all $\lambda \in
V^\wedge \setminus \{\lambda_0\}$.  By approximating $e$
we may find a polynomial $p$ with $|p(\lambda_0)| >
\frac{1}{2}$ and $|p| <\frac{1}{2}$ on $V^\wedge
\setminus\{\lambda_0\}$.  Now, necessarily, $\sup_V |p| =
|p(\lambda_0)|$. Hence the polynomial extension property
implies that there exists $g\in H^\infty({\Bbb D}^2)$ with
$\sup\limits_{{\Bbb D}^2} |g| = |g( \lambda_0)|$.  It follows
from the maximum principle that $g$ is constant.  Thus, $p$ is
constant on $V$ and $V^\wedge\setminus \{\lambda_0\}$ is
empty.  This means that, indeed, $V=\{\lambda_0\}$ and
concludes the proof of Lemma \ref{lem:5.1}.
\enddemo

Before continuing, we  remark that an alternate proof of
Lemma \ref{lem:5.1} can be obtained by showing that the
Shilov boundary of $P(V^-)$ must necessarily lie in the
boundary of ${\Bbb D}^2$ if $V$ has the polynomial
extension property and then invoking the Hugo Rossi local
maximum principle.  Either way it is curious that there
appears to be no ``elementary" proof of the lemma.

\proclaim{Lemma}\label{lem:5.2}
Let $V\subseteq {\Bbb D}^2$ and assume that $V^\wedge
\cap {\Bbb D}^2=V$.  If $V$ has the polynomial extension
property and $\lambda_0 \in {\Bbb D}^2 \setminus V${\rm ,} then
there exists $h\in H^\infty({\Bbb D}^2)$ with
$h(\lambda_0) \not= 0$ and $h|V=0$.
\endproclaim

\demo{Proof}  Fix positive $\varepsilon$ with
$\varepsilon\, <1$.  Since $\lambda_0 \in {\Bbb D}^2
\setminus V$,
$\lambda_0
\notin V^\wedge$.  Hence there exists a polynomial $p$ with
$p(\lambda_0)=1$ and 
$\sup\limits_{V} |p| \le\, \varepsilon$.
As $V$ has the poynomial extension property, there exists $G$ in
$H^\infty({\Bbb D}^2)$ with $g|_V=p$ and
$\sup\limits_{{\Bbb D}^2} |g| \le\,
\varepsilon$.  Setting $h=p-g$ we see immediately that
$h\in H^\infty({\Bbb D}^2), h(\lambda_0) \not=0$ and
$h|V=0$.  This establishes Lemma \ref{lem:5.2}
\enddemo

We recall for the reader that the intersection of an arbitrary
collection of zero sets of holomorphic functions defined on a
domain of holomorphy $D$ is in fact analytic in $D$ (see eg.
\cite{gur}).
Hence Lemma \ref{lem:5.2} has as an immediate
corollary that any relatively polynomially convex subset of
${\Bbb D}^2$ with the polynomial extension property is in
fact an analytic variety.

In the following lemma for $z\in {\Bbb C}$ and $\varepsilon\ >0$ we
let $$\Delta_z(r) = \{w\in {\Bbb C} \ |\ \ |w-z| <\, \varepsilon\}.$$  The
lemma may be regarded as an immediate consequence of the
local descriptive theory of analytic varieties and we
therefore omit its proof. We use ${\rm Zer}(h)$ to denote the zero-set
of $h$.

\proclaim{Lemma}\label{lem:5.3}
Let $h\in H^\infty({\Bbb D}^2)${\rm ,} assume that $h\not\equiv 0${\rm ,} and
fix $\lambda_0\in {\Bbb D}^2$ with $h(\lambda_0) =0$.
There exist  $n\ge 0$ and $\varepsilon_1, \varepsilon_2\ >0$ such that for all
$$z\in \Delta_{\lambda^1_0}
(\varepsilon_1)\setminus\{\lambda^1_0\}$$ there exists a
neighborhood $U$ with $U\subseteq
\Delta_{\lambda^1_0}(\varepsilon_1)\setminus\{\lambda^1_0\}$ of
$z$ and $n$ distinct holomorphic functions $f_1, \dots f_n$
defined on $U$ with the property that
$${\rm Zer}(h) \cap \left(U\times \Delta_{\lambda^2_0} (\varepsilon_2)\right)
= \bigcup\limits^n_{\ell=1} \left\{(z, f_\ell (z)) | z\in
U\right\}.$$
\endproclaim

Thus, the zero set ${\rm Zer}(h)$ of an $H^\infty$ function near a point
$\lambda_0$ consists of $n$ horizontally stretched ``sheets"
as described in Lemma \ref{lem:5.3}.  There may or may not
be a vertical sheet $\left\{\left(\lambda^1_0, z\right): z\in
{\Bbb D}\right\}$ in ${\rm Zer}(h)$, though if there is not, necessarily
$h>0$.

Armed with the preceding three lemmas we now are ready
to commence the proof of Theorem \ref{thm:1.20} from the
introduction.  Accordingly, fix a nonempty relatively
polynomially convex subset $V$ of ${\Bbb D}^2$.  If $V$ has
one of the forms\break (i)--(iv) of Theorem \ref{thm:1.20} then $V$
is a retract, i.e.\ there exists a holomorphic mapping $\rho: {\Bbb D}^2
\to {\Bbb D}^2$ with $\rho \circ \rho = \rho$ and $\ran \rho
=V$.  Clearly  if $p$ is a polynomial$, p\circ \rho$ yields a
norm-preserving extension of $p$ to ${\Bbb D}^2$ and thus
$V$ has the polynomial extension property.

To see the reverse direction assume $V$ has the polynomial
extension property.  If $V$ consists of a single point, then
$V$ is as in (i) and we are done.  Also, since by Theorem
\ref{thm:1.19},
$V$ is balanced, Proposition \ref{prop:4.10} implies that if
$V$ contains a balanced pair of points then either $B$ is a
balanced disk (in which case $V$ is as in (iii) as well as (iv)) or
$V$ is ${\Bbb D}^2$ (in which case $V$ is as in (ii)).  Thus, we
may make the following assumptions.
\begin{eqnarray}\label{eq:5.4}
 && V \quad \mbox{ contains more than one point}. \speqnu{5.4}
\\
 &&
  V \quad \mbox{ contains no pair of balanced points}.  \label{eq:5.5} \speqnu{5.5}
\end{eqnarray}

At this point  we are able to describe how the proof of
Theorem \ref{thm:1.20} will be consummated.  Say a set
$E\subseteq {\Bbb D}^2$ is an {\it extremal disk} if $E$
has one of the forms
\begin{eqnarray}\label{eq:5.6}
 && E = \{(z, f(z))\  |\  z\in {\Bbb D}\}, \speqnu{5.6}
\\
\noalign{\noindent or}
\noalign{\vskip3pt}
 \label{eq:5.7}
&&  E = \{(f(z),z)\  |\  z\in {\Bbb D}\}  \phantom{,}\speqnu{5.7}
\end{eqnarray}
for some holomorphic mapping $f:{\Bbb D}\to {\Bbb D}$.
If $E$ is an extremal disk we say {\it $E$ is type} 1  if
(\ref{eq:5.6}) holds and we say {\it $E$ is type} 2 if
(\ref{eq:5.7}) holds.  Extremal
disks are so named because they are precisely the ranges of
extremal functions for the Kobayashi extremal problems.
Also, an extremal disk is both type 1 and type 2 if and only if
it is a balanced disk.  

We shall show in Lemma \ref{lem:5.8} below, that if
(\ref{eq:5.4}) and (\ref{eq:5.5}) hold, then given any point
$\lambda \in V$, there exist extremal disks $E\subseteq V$
that are arbitrarily close to $\lambda$.  Since (\ref{eq:5.4}) implies
via Lemma \ref{lem:5.1} that $V$ has no isolated points, it will
then follow that $V$ is a union of extremal disks
(Lemma~\ref{lem:5.17}).  Finally, we shall show that $V$ cannot contain more
than one extremal disk (Lemma~\ref{lem:5.20}) and the proof of Theorem
\ref{thm:1.20} will be complete.
\advance\theoremcount by 4
\proclaim{Lemma}\label{lem:5.8}
Let $V\subseteq {\Bbb D}^2${\rm ,} assume that $V^\wedge \cap
{\Bbb D}^2=V${\rm ,} that $V$ has the polynomial extension
property{\rm ,} and that {\rm (\ref{eq:5.4})} and {\rm (\ref{eq:5.5})} hold.  If
$\varepsilon\ >0$ and $\lambda_0 \in V${\rm ,} then there exists $z_0 \in
{\Bbb D}${\rm ,} $\mu_0 \in V,$ and holomorphic $f:{\Bbb D}\to
{\Bbb D}$ with $\|\mu_0 -\lambda_0\| <\ \varepsilon$ and such that
either
\begin{eqnarray*}
{\rm (i)}\quad \{(z, f(z)) \ | \ z\in {\Bbb D}\}&\hskip-6pt \subseteq\hskip-6pt& V
\mbox{ and }\ (z_0, f(z_0)) = \mu_0 \\
\noalign{\noindent or}
\noalign{\vskip3pt}
 {\rm (ii)}\quad \{(f(z), z) \ | \ z\in {\Bbb D}\}&\hskip-6pt \subseteq\hskip-6pt& V
\mbox{ and }\ (f(z_0, z_0) = \mu_0. \end{eqnarray*}
\endproclaim

\demo{Proof} Assume that $V$ satisfies the
hypotheses of the lemma, let $\varepsilon\ >0$ and fix $\lambda_0 \in
V$.  By Lemma \ref{lem:5.2} there exists $h\in
H^\infty({\Bbb D}^2)$ with $h\not= 0$ and $V\subseteq
{\rm Zer}(h)$.  An application of Lemma \ref{lem:5.3} to $h$ yields
$n\ge 0$ and $\varepsilon_1, \varepsilon_2\ >0$ such that for all
$z_0 \in \Delta_{\lambda^1_0} (\in _1)
\setminus\{\lambda^1_0\}$ there exists a neighborhood
$U_{z_0}$ of $z_0$ with $U_{z_0}
\subseteq\Delta_{\lambda^1_0} (\in _1)
\setminus\{\lambda^1_0\}$ and $n$ distinct holomorphic
functions $f_1, \dots, f_n$ on $U_{z_0}$ with the property that
\begin{equation}\label{eq:5.9}
{\rm Zer}(h) \cap \left( U_{z_0}
\times \Delta_{\lambda^2_0} (\in _2)\right) =
\bigcup^n_{\ell=1} \{(z, f_\ell(z)) \ | \ z\in U_{z_0}\}. \speqnu{5.9}
\end{equation}

Now, by Lemma \ref{lem:5.1}, $\lambda_0$ is not an isolated
point of $V$.  Thus, either
\begin{eqnarray}
\noalign{\vskip-9pt}
&&\speqnu{5.10}\label{eq:5.10}\\
&&\mbox{there is a sequence}\ \{w_n\} \subseteq {\Bbb D} \
\mbox{with}\ w_n\to \lambda^2_0 \ \mbox{and}\
(\lambda^1_0, w_n) \in V \ \hbox{for all } n \nonumber
\\
\noalign{\noindent 
or,}
 \label{eq:5.11}
&&\quad\mbox{there is a point}\ \mu_0\in V \ \mbox{with}\ \mu^1_0
\not= \lambda^1_0 \ \mbox{and}\ \|\mu_0-\lambda_0\| <
\min\{\varepsilon, \varepsilon_1,\varepsilon_2\}.\quad \speqnu{5.11}
\end{eqnarray}
If (\ref{eq:5.10}) obtains, then it is easy to see that the
conclusion of Lemma \ref{lem:5.8} holds.  Simply note that if
$g\in H^\infty({\Bbb D}^2)$ and $g|V=0$, then
(\ref{eq:5.10}) implies that $g(\lambda^1_0, w_n) =0$ for all
$n$.  Hence since $w_n \to \lambda^2_0 \in {\Bbb D}$,
$g(\lambda^1_0, w) =0$ for all $w\in {\Bbb D}$.  But this
implies via Lemma \ref{lem:5.2} that $(\lambda^1_0, w) \in V$
for all $w\in {\Bbb D}$, and this implies that (ii) holds with
$f(z) \equiv \lambda^1_0, \mu_0 = \lambda_0$, and $z_0
=\lambda^2_0$ which establishes Lemma \ref{lem:5.8} in the
case when (\ref{eq:5.10}) holds.

To handle  the case when (\ref{eq:5.11}) holds is more subtle.
First note that since $\mu^1_0 \not= \lambda^1_0$ and
$\|\mu_0 -\lambda_0\| 
<\min \{\varepsilon, \varepsilon_1, \varepsilon_2\}$, we have
$\mu^1_0 \in \Delta_{\lambda^1_0}(\varepsilon_1)
\setminus\{\lambda^1_0\}$ and $\mu_0 \in {\rm Zer}(h) \cap
(U_{\mu^1_0}
\times \Delta_{\lambda^1_0}(\varepsilon_1))$.  Thus, \pagebreak noting that
$V\subseteq {\rm Zer}(h)$ and that Lemma~\ref{lem:5.1} implies
$\mu_0$ is not an isolated point of $V$, we see from
(\ref{eq:5.9}) (with $z_0 = \mu^1_0$), that there exist  $\ell$
and a sequence $\{z_i\}$ with $z_i\to \mu^1_0, f_\ell (z_i) \to
\mu^2_0$ such that both $z_i \in U_{\mu^1_0}$ and $(z_i,
f_\ell (z_i)) \in V$.

Let $G_0$ denote the component of
$U_{\mu^1_0}$ containing  $\mu^1_0$.
We claim that 
$(z, f_\ell(z)) \in
V$ for all $z\in G_0$.  To see this we shall use Lemma
\ref{lem:5.2}.  Thus, let $g\in H^\infty({\Bbb D}^2)$ with
$g|_V=0$.  Define a holomorphic function $\varphi$ on $G_0$ by
$\varphi(z) = g(z, f_\ell(z))$.  Indeed, $\varphi$ is a well
defined holomorphic function since (\ref{eq:5.9}) implies that
$(z, f_\ell(z)) \subseteq {\Bbb D}^2$ if $z\in U_{\mu^1_0}
\supseteq G_0$.  

Since $z_i \to \mu^1_0, z_i \in G_0$
for sufficiently large $i$.  Also, $(z_i, f_\ell(z_i)) \in V$ so that
$\varphi(z) = g(z, f(z)) =0$ for $i$ sufficiently large.  As
$\varphi$ is holomorphic and $G_0$ is connected these facts
imply that $\varphi(z) =0$ for all $z\in G_0$.  Therefore,
if $g\in H^\infty({\Bbb D}^2)$ and
$g|_V=0$, then  $g(z, f_\ell(z))=0$ whenever $z\in G_0$.  Thus,
we see from Lemma \ref{lem:5.2} that indeed $(z, f_\ell(z))
\in V$ for all $z\in V$.

Recapping, we have shown that there is a $\mu_0 \in V$,
a connected neighborhood $G_0$ of $\mu^1_0$, and a
holomorphic $f_\ell : G_0\to {\Bbb D}$ so that
$\|\mu_0 - \lambda_0\| <\varepsilon$, 
$$\{(z, f(z))\ | \ z\in G_0\} \ \subseteq \ V $$
and $$
\Big(\mu^1_0, f_\ell (\mu^1_0)\Big) \ =\ \mu_0.$$
Notice how similar this is to (i) in the conclusion of Lemma
\ref{lem:5.8} (with $z_0 = \mu^1_0$).  Indeed, we would have
established the lemma if it were the case that $f_\ell$ had an
analytic continuation $\tilde{f}_\ell$ to the entire disk with
$\ran \tilde{f} \subseteq {\Bbb D}$.  This however does not
have to be the case in general.
The next part of the proof of Lemma \ref{lem:5.8} addresses
this issue.

Define a set $S$ in ${\Bbb C}^2$ by letting $S=\{(z, w) \ |\ z,
w\in G_0 $ and $z\not= w\}$ and define a real-valued function
$\rho$ on $S$ by  letting 
$$\rho (z, w) \ =\ d(z, w) - d(f_\ell (z),
f_\ell(w)).
$$
Since $G_0$ is connected and open, $S$
is also
connected.  Furthermore, notice that since $(z, f_\ell(z)) \in
V$ whenever $z\in G_0$, (\ref{eq:5.5}) implies that $\rho(z,
w) \not= 0$ for all $(z, w) \in S$.  Since $\rho$ is continuous it
follows from the intermediate value theorem that either
\begin{equation}\label{eq:5.12}
\rho(z, w) >0 \quad \mbox{for all}\quad (z, w) \in S, \speqnu{5.12}
\end{equation}
or
\begin{equation}\label{eq:5.13}
\rho(z, w) < 0 \quad \mbox{for all}\quad (z, w) \in S. \speqnu{5.13}
\end{equation}

Now, if (\ref{eq:5.12}) holds then, indeed, $f_\ell$ does
extend to a holomorphic mapping $f:{\Bbb D} \to
{\Bbb D}$.  To see this consider function elements $(G,
f_G)$, with $G$ a connected open subset of ${\Bbb D}$ such
that $G_0 \subseteq G$ and $f_G: G\to {\Bbb D}$ a
holomorphic function with $f_G\ |_{G_0} =f_\ell$.

We first claim that
\begin{equation}\label{eq:5.14}
(z, f_G(z)) \in V\quad \mbox{for all}\quad z\in G. \speqnu{5.14}
\end{equation}
To see (\ref{eq:5.14}), let $G_1$ be the component of $\{z\in
G \ | \ (z, f_G(z)) \in V\}$ that contains $G_0$.  The proof of
(\ref{eq:5.14}) will follow if we can show that it is
both open and closed as a subset of $G$.  That $G_1$ is open
follows from Lemmas \ref{lem:5.2} and \ref{lem:5.3}.  That
$G_1$ is closed follows from the fact that $V^\wedge
\cap{\Bbb D}^2=V$.

Next we claim that if $ (G, f_G)$ is a function element then
\begin{equation}\label{eq:5.15}
f_G\ \mbox{is strictly contractive on}\ G\ \mbox{in the pseudo-hyperbolic metric}. \hskip.25in  \speqnu{5.15}
\end{equation}
To prove (\ref{eq:5.15}) note that we are assuming that
(\ref{eq:5.12}) holds, i.e., $f_G$ is contractive on $G_0$.  Just
as in the proof that either (\ref{eq:5.12}) or (\ref{eq:5.13})
hold we see that $\rho >0$ on  the set $\{(z, w) \ |\ (z,w) \in G$ and
$z\not= w\}$.  Hence (\ref{eq:5.15}) obtains.
\vglue4pt

We next claim that if $(G_1, f_{G_1}), (G_2, f_{G_2})$ are two
function elements, then
\begin{equation}\label{eq:5.16}
f_{G_1} = f_{G_2} \quad \mbox{on}\quad G_1 \cap G_2. \speqnu{5.16}
\end{equation}
To see (\ref{eq:5.16}) we argue by contradiction.  Thus,
assume that $z_1 \in G_1 \cap G_2$ with $f_{G_1}(z_1) \not=
f_{G_2} (z_1)$.  Choose a curve $z:[0, 1] \to G_2$ with $z(0) \in
G_0$ and $z(1) =z_1$.  Define $\psi: [0, 1] \to {\Bbb R}$ by
$$
\psi(t) \ = \ d(z(t), z_1) - d(f_{G_2} z(t), f_{G_1}(z_1)).
$$
When
$t=0, f_{G_2} (z_t) = f_{G_2}(z_0) = f_{G_1}(z_0)$ since
$f_{G_1} = f_{G_2}$ on $G_0$ and $z_0 \in G_0$.  Hence
$\psi(0) >0$ by (\ref{eq:5.15}).  On the other hand, $\psi(1) <0$
since $z(1) =z_1 $ and $f_{G_1}(z_1) \not= f_{G_2}(z_1)$.  Hence
there exists $t_0$ such that $\psi(t_0)=0$.  But this says that
the pairs of points $(z(t_0), f_{G_2}(z(t_0)))$ and $(z_1,
f_{G_1}(z_1))$, which by (\ref{eq:5.14}) are in $V$, are balanced
pairs.
This contradiction to (\ref{eq:5.5}) establishes (\ref{eq:5.16}).

Finally, note that (\ref{eq:5.16}) allows us to construct a
maximal function element $(G, f_G)$.  We claim that
$G={\Bbb D}$.  
Indeed, if $w\in \partial G \cap {\Bbb D}$, one of two things
must have happened. Either
$$
\lim_{ G \ni z \to w} \, |f(z)| = 1 ,
$$
which would contradict \ref{eq:5.15}; or 
$f_G$ has a singularity at $w$.
By Lemmas \ref{lem:5.2} and
\ref{lem:5.3}, this singularity is isolated.
Since isolated singularities of bounded holomorphic
functions are removable, $\partial G \cap
{\Bbb D}=\emptyset$.

We are finally able to complete the proof of Lemma
\ref{lem:5.8}.  If (\ref{eq:5.12}) holds, then take $f=f_G$ where
 $G$ is as in the previous paragraph. Now, (i) holds by
(\ref{eq:5.14}).  If (\ref{eq:5.13}) holds then execute the
argument just carried out under the assumption that
(\ref{eq:5.12}) held but with the coordinates interchanged,  {i.e.}\
$G_0$ replaced with an appropriate connected neighborhood
of $\mu^2_0, f_{G_0}$ replaced with $f^{-1}_\ell$ and the
manifolds $(z, f_G(z))$ replaced by $(f_G(z), z)$.  This yields the fact
that (ii) holds and completes the proof of Lemma
\ref{lem:5.8}.
\enddemo
\advance\theoremcount by 8
\proclaim{Lemma}\label{lem:5.17}
If $V$ satisfies the hypotheses of Lemma {\rm \ref{lem:5.8},} then
$V$ is a union of extremal disks.
\endproclaim  

\demo{Proof}  Fix $\lambda \in V$.  By Lemma
\ref{lem:5.8},  there exist  sequences $\{z_n\}, \{\lambda_n\}$,
and $\{f_n\}$ such that $\lambda_n \to \lambda$ and such that
either
\begin{equation}\label{eq:5.18}
\{(z, f_n(z)) | z\in {\Bbb D}\} \subseteq V\quad
\mbox{and}\quad (z_n, f_n(z_n)) = \lambda_n\quad
\mbox{for all}\ n, \hskip.5in \speqnu{5.18}
\end{equation}
or
\begin{equation}\label{eq:5.19}
\{(f_n(z), z) | z\in {\Bbb D}\} \subseteq V\quad
\mbox{and}\quad (f_n(z_n), z_n) = \lambda_n\quad
\mbox{for all}\ n. \hskip.4in \speqnu{5.19}
\end{equation}
Suppose (\ref{eq:5.18}) holds.  Evidently, $z_n\to \lambda^1$
and $f_n(z_n) \to \lambda^2$.  Furthermore, since $\sup |f_n|
\le 1$, by passing to a subsequence if necessary we may
assume that $f_n\to f$ uniformly on compact subsets of
${\Bbb D}$.  Necessarily, $$
\{(z, f(z)) \ | \ z\in
{\Bbb D}\}\  \subseteq  \ V
$$
(since $V$ is a closed subset of
${\Bbb D}^2$) and obviously, $(\lambda^1, f(\lambda^1))
=\lambda$.  Thus, $\lambda$ is in a subset of $V$ which is an
extremal disk of type 1.  Likewise, if we assume that
(\ref{eq:5.19}) holds, then $\lambda$ is in a subset of $V$
which is an extremal disk of type~2.  This completes the proof
of Lemma \ref{lem:5.17}.
\enddemo

In light of Lemma \ref{lem:5.17} the proof of Theorem
\ref{thm:1.20} will be complete once we have established our
final lemma.
\advance\theoremcount by 2

\proclaim{Lemma}\label{lem:5.20}
If $V$ satisfies the hypotheses of Lemma {\rm \ref{lem:5.8},} then
$V$ cannot contain more than  one extremal disk.
\endproclaim

\demo{Proof}  There are two cases to rule out: $V$
contains two extremal disks of the same type and $V$
contains two extremal disks of different type.

First suppose $V$ contain two distinct disks of type 1:
$$E_f = \{(z, f(z)) | z\in {\Bbb D}\} \ \mbox{and}\ E_g = \{(z,
g(z)) |z\in {\Bbb D}\}.$$
We claim that
\begin{equation}\label{eq:5.21}
d(z, w) < d(f(z), g(w)) \ \mbox{whenever}\ z\not= w\ \mbox{or}\
f(z)
\not= g(w).  \hskip.5in\speqnu{5.21}
\end{equation}
To prove (\ref{eq:5.21}) let $S=\{(z, w)| z\not= w $ or $f(z)
\not= g(w)\}$, define $\rho: S\to {\Bbb R}$ by $\rho(z, w)
=d(f(z), g(w)) -d(z, w)$.  Since $E_f \not= E_g, f \not= g$ and
we see that for some $z_0, (z_0, z_0) \in S$ and $\rho(z_0,
z_0) >0$.  Also, by (\ref{eq:5.5}), $\rho(z, w) \not= 0$ for all
$(z, w) \in S$.  It follows by the fact that $S$ is connected
and $\rho$ is continuous that $\rho >0$ on all of $S$, i.e.,
(\ref{eq:5.21}) holds.

Armed with (\ref{eq:5.21}) we are able to show that $f=g$ in
the following way.  Alternately fixing $w$ and letting $z\to
\partial {\Bbb D}$ and fixing $z$ and letting $w\to \partial
{\Bbb D}$ immediately imply via (\ref{eq:5.21}) that $f$
and $g$ are inner.  Furthermore since $|f|$ and $|g|$ have
unimodular limits along all paths to $\partial {\Bbb D}$, in
fact $f$ and $g$ are finite Blaschke products.  Observe that if
$z$ is a zero of $f$ and $w$ is a zero of $g$, then (\ref{eq:5.21})
implies that $z=w$.  Thus, $f$ and  $g$ are automorphisms of
${\Bbb D}$.  Letting $w_n \to g^{-1} (f(z))$ in such a way
that $w_n \not= z$ and $w_n \not= g^{-1} (f(z))$ reveals via
(\ref{eq:5.21}) that $g^{-1} (f(z)) =z$, i.e., $f=g$.  Thus, we see
that $E_f$ and $E_g$ are not distinct after all, a contradiction
that establishes the fact that $U$ does not contain two extremal disks of
type 1.  Of course a nearly identical argument shows that $V$
cannot contain two extremal disks of type 2 as well.

Now suppose $V$ contains two distinct extremal disks of
opposite type say $E_f = \{(z, f(z)) |z\in {\Bbb D}\}$ and
$E_g = \{(g(z), z) | z\in {\Bbb D}\}$.  The argument of the
previous part of the proof gives that either \advance\eqcount by 21
\begin{equation}\label{eq:5.22}
d(z, g(w)) < d(f(z), w) \ \mbox{whenever}\ z \not= g(w) \
\mbox{or}\ w\not= f(z) \hskip.6in
\end{equation}
or
\begin{equation}\label{eq:5.23}
d(f(z) w) < d(z, g(w)) \ \mbox{whenever}\ z \not= g(w) \
\mbox{or}\ w\not= f(z). \hskip.6in\speqnu{5.23}
\end{equation}
If (\ref{eq:5.22}) holds, then choosing $w_n \to f(z)$ in such a
way that $w_n \not= f(z)$ and $g(w_n) \not= z$ gives that
$g(f(z)) =z$.   But this implies that $E_f =E_g$, a
contradiction.  Likewise, if (\ref{eq:5.23}) holds, we reach a
contradiction  and Lemma
\ref{lem:5.20} is established.
\enddemo 
 
\section{Appendix by Pascal J. Thomas}

The following direct function theoretic proof of Theorem~\ref{thm:1.19} is
due to Pascal J. Thomas, of Universit\'e Paul Sabatier.
It is based in part on an argument, in the context of the unit
ball,
that Jean-Pierre Rosay showed Thomas in a conversation in 1993.

\demo{Proof} 
First, note that $(p_1,p_2)$ is a balanced pair of points
if and only if
there
exist $m_1, m_2 \in \cal M$, two automorphisms of the unit disk such
that
$m_1 (p_1^1) = m_2 (p_2^1) =0$ and $m_1 (p_1^2) = m_2 (p_2^2)$.
Polynomials themselves are not invariant under automorphisms of the
bidisk, but any composition of a polynomial by such an automorphism
is
easily seen to be uniformly approximable by polynomials on the
closed
bidisk. So we
might reduce ourselves to the case $p_1=(0,0), p_2=(\lambda,
\lambda)$.
In that
case, the unique totally geodesic complex disk going through both
points
is the diagonal $\Delta := \{ (z,z) : z \in \Bbb D \}$.

By the hypothesis that $\overline V$ is polynomially convex, to prove
that
$\Delta \subset V$, it will be enough to prove that $\overline V
\supset \partial
\Delta := \{ (z,z) : z \in \partial \Bbb D \}$. Suppose not; then
there
exists $\theta_0$ and
$\delta > 0 $ such that $\overline V \subset \overline{\Bbb D^2}
\setminus
D(e^{i\theta_0},\delta)^2 =:K$.

Consider the holomorphic retraction from $\Bbb D^2$ onto $\Delta$
(or rather, the unit disk parametrizing $\Delta$) given by
$\pi (z,w) := \frac12 (z+w)$. By strict convexity of the unit disk,
$\pi^{-1} \{e^{i\theta_0}\} \cap \overline{\Bbb D^2}
= \{(e^{i\theta_0},e^{i\theta_0})\}$, so that
$(e^{i\theta_0},e^{i\theta_0})
\notin \pi(K)$ and by compactness, there exists $\epsilon > 0$ so
that
$\pi(V) \subset \Bbb D \setminus D(e^{i\theta_0},\epsilon)
=: \Omega$.

Now there exists a holomorphic function $f$ on $\Omega$
such that $f(\Omega) \subset \Bbb D$, $f(0)=0$ and $|f(\lambda)| >
|\lambda|$; take for instance $f$ the conformal mapping from
$\Omega$ to
$\Bbb D$, the property is then simply the Schwarz lemma applied to
the
map $f^{-1}$ at the point $f(\lambda)$. I claim that the map $(f
\circ
\pi)|_V$, which has uniform norm bounded by $1$, cannot extend to
any
$\tilde f \in H^\infty (\Bbb D^2)$ with $\|\tilde f\|_\infty \le
1$.

Indeed, if such an $\tilde f$ existed, then $f_1(z) := \tilde
f(z,z)$
would give a function from the unit disk to itself such that
$f_1(0)=0$,
but by the extension property $f_1(\lambda) = f(\lambda)$, which
violates
the Schwarz lemma.

To see that the {\it polynomial} extension property is violated, it
is
enough
to approximate $f$ by a polynomial $P$ (by Mergelyan's theorem), so
that\break
$\sup_{z \in \Omega} |P(z)| \le 1$ and $|P(\lambda)| > |\lambda|$,
and
then run the same argument with $(P \circ \pi)|_V$.
\enddemo

 \end{document}